\documentclass[11pt,a4paper, notitlepage]{article}
\usepackage[a4paper] {geometry}
\usepackage{amsmath,amsfonts,amssymb,mathrsfs,amsthm}
\usepackage{multirow}
\usepackage[english]{babel}
\usepackage[all]{xy}
\usepackage{tikz}
\numberwithin{equation}{section}
\numberwithin{figure}{section}
\newtheorem{thm}{Theorem}[section]

\theoremstyle{definition}
\newtheorem{defx}[thm]{Definition}
\newcommand{\C}{I\!\!\!C}
\newcommand{\h}{I\!\!H}

\newcommand{\R}{I\!\!R}
\newcommand{\N}{I\!\!\!N}

\newcommand{\be}{\begin{enumerate}}
\newcommand{\ee}{\end{enumerate}}
\newcommand{\bq}{\begin{eqnarray*}}
\newcommand{\eq}{\end{eqnarray*}}

\DeclareMathSizes{12}{13}{10}{10}

\title{Invariant Differential Operators on the Heisenberg Group}
\author{M.E. Egwe\\
\small Department of Mathematics,\\
\small University of Ibadan,\\
\small Ibadan,  Nigeria, $200005$\\
$murphy.egwe@ui.edu.ng$}
\begin{document}
\maketitle

\begin{abstract}
In this paper, we elucidate certain properties of the $(2n+1)$-dimensional Heisenberg group, and establish some theorems on the invariant differential operators on the group.
\end{abstract}

{\bf Keywords:} $Heisenberg \;group,\;Invariant\;Differential Operator,\; Heisenberg\; Laplacian\\Factorization,\; Universal\; enveloping\; algebra,\; solvability.$

\section{Preliminaries}
The Heisenberg group (of order $n$), $\h_{n}$ is a
noncommutative nilpotent Lie group whose underlying manifold is $\C^{n}\times \mathbb{\R}$
with coordinates $(z,t)=(z_{1},z_{2},...,z_{n},t)$ and group law given by
$$ (z,t)(z',t')=(z+z',t+t'+2\Im m\; z.z')\;\; \mbox{where}\;\;
 z.z'=\sum_{j=1}^{n}z_{j}\bar z_{j}\;\;\;\;z\in\C^{n},\;\;t\in\mathbb{\R}. $$
Setting $z_{j}=x_{j}+y_{j}$, then $(x_{1},x_{2},...,x_{n},y_1,y_2,\cdots,y_n,t)$ forms a real coordinate
system for $\h_{n}$.
In this coordinate system, we define the following vector fields:
$$ X_{j}=\frac{\partial}{\partial x_{j}}+2y_{j}\frac{\partial}{\partial t},\;\; Y_{j}=\frac
{\partial}{\partial y_{j}}-2x_{j}\frac{\partial}{\partial t},\;\;T=\frac{\partial}{\partial t}. $$ It is clear from \cite{FolSte} that $\{X_{1},X_{2},...,X_{n},Y_{1},Y_{2},...,Y_{n},T\}$ is a basis for the left invariant vector fields on $\h_{n}.$ These vector fields span the Lie algebra $\mathfrak{h}_n$ of $\h_n$ and
the following commutation relations hold:
 $$[Y_{j},X_{k}]=4\delta_{jk}T,\;\;\; [Y_{j},Y_{k}]=[X_j,X_k]=[X_{j},T]=[Y_{j},T]=0.$$
Similarly, we obtain the complex vector fields by setting
$$\left.\begin{array}{lr}
Z_j =\frac{1}{2}(X_j-iY_j) = \frac{\partial}{\partial z_j} +
i\bar{z}\frac{\partial}{\partial t}\\
\ \\
\bar{Z}_j =\frac{1}{2}(X_j+iY_j) = \frac{\partial}{\partial
\bar{z}_j} -
iz\frac{\partial}{\partial t}\end{array}\right\}.\eqno(1)$$
In the complex coordinate, we also have the commutation relations
$$[Z_j,\bar{Z}_k] = -2\delta_{jk}T,\;\;[Z_j,Z_j] = [\bar{Z}_k,\bar{Z}_k]= [Z_j,T]=[\bar{Z},T] = 0.$$
The Haar measure on $\h_n$ is the Lebesgue measure
$dzd\bar{z}dt$ on $\C^n\times\mathbb{\R}$ \cite{Howe}. In particular, for $n=1$, we obtain the $3$-dimensional Heisenberg group $\h_1\cong\mathbb{\R}^3$ (since $\C^n\cong\mathbb{\R}^{2n})$. Hence $\h_n$ may also be referred to as $(2n+1)$-dimensional Heisenberg group.

For a real non-zero $\lambda$, let $\cal H^{\lambda}$
be the Hilbert space of entire functions $F$ on $\C^{r}$ such that
$$\int_{\C^{r}}e^{-2|\lambda||z|^{2}}|F(z)|^{2}dz=|\lambda|^{-r}\| F\|^2_{\cal H^{\lambda}}$$
is finite. The monomials $z^{n}=z_{1}^{n_{1}}...z_{r}^{n_{r}}, (n_{j}\in \N,$ and $n=n_1+n_2+\cdots n_r$) are orthogonal
in $\cal H^{\lambda}$ and the functions
$$\phi_{n}^{\lambda}(z)=(n!)^{-1/2}(\frac{2}{\pi})^{r/2}(2|\lambda|z)^{n}\;\;\; (n!=n_1!,...,n_r!)
$$ form an orthonormal basis for $\cal H^{\lambda}$.\\
\\
The group $\h_{n}$  acts on $\cal H^{\lambda}$ for positive $\lambda$  by
$$U_{(z,t)}^{\lambda}F(w)=e^{-i\lambda t+\lambda (2{\langle w,z\rangle -|z|^2})}F(w-z).$$
and for negative $\lambda$ by
$$U_{(z,t)}^{\lambda}F(w)=e^{-i\lambda t-\lambda (2{\langle w,\bar z\rangle -|z|^2})}F(w-\bar z).$$
For $\lambda =0$ we get the one dimensional representations
$$\chi_{w}(z,t)=e^{-iRe\langle w,z\rangle}$$
These are the irreducible unitary representations of $\h_n$ up to equivalence.

One significant structure that accompanies the Heisenberg group is the family of dilations
$$\delta_{\pm \lambda}(z,t)=(\pm \lambda z,\pm \lambda^2t),\;\; \lambda >0 $$
This family is an automorphism of $\h_n.$ Now, if $\sigma:\C\rightarrow \C$ is an automorphism, there exists an induced automorphism, $\widetilde{\sigma} \in Aut(\h_n),$ such that
$$\widetilde{\sigma}(z,t)=(\sigma z, t).$$ For simplicity, assume that $\widetilde{\sigma}$ and $\sigma$ coincide. Thus we may simply assume that if $\sigma\in Aut(\h_n),$ we have $\sigma (z,t)=(\sigma z, t).$\\
\section{Differential Operators}
\begin{defx}
A linear partial differential operator with smooth coefficient (of order m) is an operator of the form
$$Lf = \sum_{|\alpha|=m}a_{\alpha}\partial^{\alpha}f,\;\;\;a_{\alpha}\in \mathbb{C}(\mathbb{\R}^n).$$
If $a_{\alpha}'s$ are constants, we call $L$ a constant - coefficient operator.

Next, let $G$ be a real or complex Lie group and $K=\mathbb{\R}$ or $\C$.  Let $C^{\infty}(G;K)$ denote the K-valued $C^{\infty}$-functions on $G$.  A left invariant differential operator on $G$ is a linear map
$$D:C^{\infty}(G;K)\longrightarrow C^{\infty}(G;K),$$
so that
\be
\item[(i)] if $g\in G$, then $D\circ L _g=L_g\circ D,\;\;(L_gf(x))=f(g^{-1}x),\;\; f\in C^{\infty}(G;K);$
\item[(ii)] if $g\in G$, then there is a neighbourhood $U$ of $g$ and a system of coordinates $\{x_1,\cdots,x_n\}$ (holomorphic if $K=\C)$ so that
$$Df|_{_U}=\sum_{|i|=m}a_{i_1\cdots i_n}\frac{\partial^{i_1+\dots+i_n}}{\partial x^{i_1}_1+\cdots + \partial x^{i_n}_n}f,$$
$f\in C^{\infty}(G;K)$, where $a_{i_1\cdots i_n}\in C^{\infty}(U;K)$ and depends only on $D$ and $U$.\ee

An operator that occurs as an analogue (for the Heisenberg group) of the Laplacian
$\Delta=\displaystyle\sum_{j=1}^n\frac{\partial^2}{\partial {x_j}}$ on $\mathbb{\R}^n$ is denoted by ${\cal L}_\alpha$
where $\alpha$ is a parameter and defined by
$${\cal L}_\alpha=-\frac{1}{2}\displaystyle\sum_{j=1}^n(Z_j\bar{Z}_j+\bar{Z}_jZ_j)+i\alpha T, \eqno (3.1)$$
where
$${Z}_j=\frac{1}{2}(X_j+Y_j)\;\mbox{and}\; \bar{Z}_j=\frac{1}{2}(X_j-Y_j),$$
so that ${\cal L}_\alpha$ can be written as
$${\cal L}_\alpha=\frac{1}{4}\displaystyle\sum_{j=1}^n(X_{j}^2+Y_{j}^2)+i\alpha T.$$
${\cal L}_\alpha$ is called the \emph{sublaplacian.}

${\cal L}_\alpha$ satisfies symmetry properties analogous to those of $\Delta$ on $\mathbb{\R}^n.$ It is left-invariant on $\h_n,$ of degree $2$ with respect to the dilation automorphism of $\h_n$ and is invariant under unitary rotations(\cite{Egwe1},\cite{Egwe2},\cite{Egwe3},\cite{Egwe4},\cite{Egwe5},\cite{Egwe6},\cite{Egwe7},\cite{Egwe8},\cite{Egwe9}).\\
\end{defx}
A differential operator $L$ on $\mathbb{\R}^n$ is said to be locally solvable at $x\in \mathbb{\R}^n$ if for every $f\in\C_0^\infty(\mathbb{\R}^n),$ there is a function (or a distribution)  $u$ and a neighbouhood $U$ of $x$ such that $Lu=f$ holds (in the sense of distribution) on U. Not all partial differential operators are solvable. For example, in \cite{Egwe10}, it was shown that the Heisenberg Laplacian (i.e., ${\cal L}_\alpha,\;\; \alpha=0$) of the form
\begin{eqnarray*}
\Delta_{\h_n}&:=& \sum^n_{j=1}X_j\circ X_j+Y_j\circ Y_j\\
&=&\sum^n_{j=1}\frac{\partial^2}{\partial
x^2_j}+\frac{\partial^2}{\partial
y^2_j}+4y_j\frac{\partial^2}{\partial x_j\partial
t}-4x_j\frac{\partial^2}{\partial y_j\partial
t}+4(x^2_j+y^2_j)\frac{\partial^2}{\partial t^2}
\hspace{0.5in} (3.2) \end{eqnarray*}
can be expressed in the factorized form
$AA^\dag$
where
$$A =\left(\displaystyle\frac{\partial}{\partial x}-i \frac{\partial}{\partial
y}+(2y-2ix)\frac{\partial}{\partial
t}\right) \mbox{and}\; A^{\dag} =\left(\displaystyle\frac{\partial}{\partial
x}+i\frac{\partial}{\partial y}+(2y+2ix)\frac{\partial}{\partial
t}\right).\eqno(3.3)$$ Consequently, it was shown that this operator is nowhere solvable and a trivial group-invariant solution was constructed on the group.

In this paper, we shall give the necessary and sufficient condition for local solvability of the factors in (3.3).

Conditions for local solvability of partial differential operators are well-known in literature. For example,\cite{Muller1},\cite{Muller2},\cite{Winfield}, \cite{Hormander},\cite{Rothschild1}, \cite{Rothschild2},\cite{Egwe10} and the references therein.
Let $\lambda\in\mathbb{\R}\backslash\{0\}=\mathbb{\R}^*$.  Consider the representations $\pi_{\lambda}$ of $\h_n$ on $L^2(\mathbb{\R}^n)$
$$[\pi_{\lambda}(x,y,t)\phi](u)=e^{i\lambda(t+(u+\frac{x}{2})y)}\phi(u+x),\phi\in L^2(\mathbb{\R}^n).$$
If $f$ is an integrable function on $\h_n$, then the integrated representation is given by
$$\begin{array}{rcl}
[\pi_{\lambda}(f)\phi](u)&=&\displaystyle\int_{_{\h_n}}f(x,y,t)e^{i\lambda(t+(u+\frac{x}{2})y)}\phi(u+x)dxdydt\\
\ \\
&=&\displaystyle\int_{_{\mathbb{\R}^n}}\bar{{\cal F}}_{y,t}f(x,\lambda(u+\frac{x}{2}),\lambda)\phi(x+u)dx\\
\ \\
&=&\displaystyle\int_{_{\mathbb{\R}^n}}\bar{{\cal F}}_{y,t}f(v-u,\lambda(v+u)/2,\lambda)\phi(u)du\\
\ \\
&=&\displaystyle\int_{_{\mathbb{\R}^n}}K^{\lambda}_f(v,u)\phi(u)du;\end{array}$$
here $\bar{{\cal F}}_{y,t}f$ denotes the partial Euclidean Fourier co-transform in the indicated variables.

Now recall that the irreducible unitary representations of $\h_n$ are given by:
\be
\item[(i)] Finite-dimensional: $V=\C$, and the representations are parameterized by $(a,b)\in\mathbb{\R}^{2n}$ with
$$\pi_{(a,b)}(X_j)=ia_j,\;\;\pi_{(a,b)}(Y_j)=ib_j,\;\mbox{and}\;\pi(T)=0.$$
\item[(ii)] Infinite-dimensional: $V=L^2(\mathbb{\R}^n)$, and the representations are parameterized by $\lambda\in\mathbb{\R}^*$ with \ee
$$\begin{array}{rcl}
\pi_{\lambda}(X_j)&=&\displaystyle|\lambda|^{1/2}\frac{\partial}{\partial x_j}, \pi_{\lambda}(Y_j)=i|\lambda|^{1/2}(sgn\;\lambda)x_j,\;\;
\mbox{and}\;\pi_{\lambda}(T)=i\lambda.\end{array}\eqno (3.7)$$
Here, $X,Y \mbox{and}\;\; T$ are vector fields and $j=1,2,\cdots,n.$\\
\begin{thm}
Let $\pi_\lambda$ be a representation of ${\cal U}(\mathfrak{H}_n)$ in $L^2(\mathbb{\R}^n),$
and let $P$ be a homogeneous left invariant differential operator on the Heisenberg group $\h_n$.  Then $P$ is locally solvable if $\pi_\lambda(P)$ has a bounded left inverse for every infinite dimensional irreducible unitary representation $\pi$ of $\h_n$.
\end{thm}
\textbf{Proof:} It suffices to prove that there exists $N>0$ such that for every $v\in C^{\infty}_c(\h_n)$
$$Pu=T^N_v,\eqno (3.8)$$
has a solution $u\in C^{\infty}(\h_n).$\\

Since $T=\displaystyle\frac{\partial}{\partial t}$ is in suitable coordinates, one may obtain $\nu_1\in C^{\infty}(\h_n)$ satisfying
$$T^N_{\nu_1}=f$$
by integrating $f\;N$ times in $t$.  Now if $U\subset C^{\infty}(\h_n)$ has compact closure, choose $\nu\in C^{\infty}_0(\h_n)$ such that $\nu\equiv\nu_1$ on $U$.  Then if $u\in C^{\infty}(\h_n)$ satisfies (3.7),
$$Pu=T^N_{\nu_1}=f\in U.$$

We now prove the theorem by constructing a solution of (3.8), using the Plancherel theorem for $\h_n$.  If $\varphi\in C^{\infty}_0(\h_n)$, and $\pi_{\lambda}$ is an infinite dimensional irreducible unitary representation of $\h_n$, then the operator
$$\pi_{\lambda}(\varphi)=\int_{\h_n}\varphi(x)\pi_{\lambda}(x)dx$$
is of trace class.  In particular, if $\{e_i\}$ is any orthonormal basis of $L^2(\mathbb{\R}^n)$, then
$$tr(\pi_{\lambda}(\varphi))=\sum_i(\pi_{\lambda}(\varphi)e_i,e_i)<\infty.$$
Now suppose $P$ is homogeneous of degree $\alpha>0$.  Then
$$\pi_{\lambda}(P)=\left\{\begin{array}{l}
|\lambda|^{\alpha/2}\pi_\lambda(P_*), \;\;\lambda>0\\
\ \\
|\lambda|^{\alpha/2}{\pi_{\lambda}}^{-1}(P), \;\;\lambda<0\end{array}\right.\eqno (3.9)$$
by (3.8).
The hypothesis of the theorem is then that $\pi_\lambda(P)$ and ${\pi_\lambda}^{-1}(P)$ have bounded inverses, to be denoted $B_1$ and $B_{-1}$.  Following [102], we first reduce to the case where $B_1$ and $B_{-1}$ are of trace class.  This is accomplished by replacing $P$ by $PQ$, where
$$Q=(\displaystyle\sum_j X^2_j+Y^2_j)^{N'}$$
for large $N'$.  Indeed from [102, Lemma 3.3], we have that $\pi_1(Q)^{-1}$ is of trace class.  Since the product of a bounded operator and an operator of trace class is again of trace class, the inverse of $\pi_1(PQ)=\pi_1(P)\pi_1(Q)$ is of trace class.  $PQ$ again satisfies the hypothesis of the theorem and the local solvability of $PQ$ implies that of $P$.  Replacing $P$ by $PQ$ if necessary, we may assume that $B_1$ and $B_{-1}$ are of trace class.  To solve (3.8), let $u$ be the linear mapping
$$u:C^{\infty}_0(\h_n)\longrightarrow\C$$
defined  by
$$\varphi\longrightarrow u(\varphi)=\int_{_{\mathbb{\R}^*}}(i\lambda)^Ntr(\pi_{\lambda}(\tilde{\nu})\pi_{\lambda}
(\varphi))(\pi_{\lambda}(P))^{-1})|\lambda|^nd\lambda\eqno (3.10)$$
for some large $N$, where $\tilde{\nu}(x)=\nu(x^{-1})$.  We shall first show that $u$ is an $L^2-$function by proving
$$|u(\varphi)|\leq C\|\varphi\|_{L^2}.\eqno (3.11)$$
By (3.9), the absolute value of the right-hand side of (3.10) is bounded by\\
$\displaystyle|\int_{\lambda>0}|\lambda|^{N-\alpha/2}|(tr(\pi_{\lambda}(\tilde{\nu})\pi_{\lambda}(\varphi)
B_1))||\lambda|^nd\lambda|$
$$+|\int_{\lambda<0}|\lambda|^{N-\alpha/2}|tr(\pi_{\lambda}(\tilde{\nu})\pi_{\lambda}(\varphi)B_{-1})
||\lambda|^n |d\lambda|.\eqno (3.12)$$
If $\alpha>0$ is arbitrary, then\\
$\displaystyle|\int_{\lambda>0}|\lambda|^{N-\alpha/2}|(tr(\pi_{\lambda}(\tilde{\nu})\pi_{\lambda}(\varphi)
B_1))||\lambda|^nd\lambda|$
$$\leq\int_{0<\lambda<\alpha}|\lambda|^{N-\alpha/2}|tr(\pi_{\lambda}(\tilde{\nu})\pi_{\lambda}(\varphi)B_1))|
|\lambda|^n|d\lambda|$$
$$+|\int_{\alpha\leq\lambda<\infty}|\lambda|^{N-\alpha/2}|tr(\pi_{\lambda}(\tilde{\nu})\pi_{\lambda}
(\varphi)B_1))|\lambda|^n|d\lambda|.\eqno (3.13)$$

The generalized Schwartz inequality for arbitrary operators $C,D$ of trace class states that
$$|tr(CD)|\leq (tr(CC^*))^{\frac{1}{2}}(tr(DD^*))^{\frac{1}{2}}.\eqno (3.14)$$
The first integral on the right of (3.13) is bounded by
$$\left[\int^a_0|\lambda|^{2(N-\alpha/2)+n}d\lambda\right]^{\frac{1}{2}}\left[\int^{\lambda}_0|tr(\pi_{\lambda}
(\tilde{\nu})\pi_{\lambda}(\varphi)B_1)|^2|\lambda|^nd\lambda\right]^{\frac{1}{2}}\eqno (3.15)$$
by the Standard Schwartz inequality.  If $2(N-\alpha/2)+n>-1$ i.e., if
$$N>(-(n+1)+\alpha)/2,$$
then the first factor in (3.15) is finite, say bounded by $C>0$.  By (3.14), (3.15) is bounded by
$$C\left[tr(B_1B^*_1)\right]^{\frac{1}{2}}\left[\int^{\lambda}_0tr(\pi_{\lambda}(\tilde{\nu})\pi_{\lambda}
(\varphi)\pi_{\lambda}(\varphi^*)\pi_{\lambda}(\tilde{\nu})^*)|\lambda|^nd\lambda\right]^{\frac{1}{2}}.\eqno (3.16)$$

Now recall that the group convolution of, say, two functions $\psi, \chi\in C^{\infty}_0$ is defined as
$$(\psi*\chi)(g)=\int\psi(x^{-1})\chi(xg)dx.\eqno (3.17)$$
Recall from 1.4.6 that
$$\pi_{\lambda}(\psi*\chi)=\pi_{\lambda}(\psi)\pi_{\lambda}(\chi).$$
Hence (3.16) is bounded by
$$C'\left|\int^{\lambda}_0tr(\pi_{\lambda}(\tilde{\nu}*\varphi)\pi_{\lambda}(\tilde{\nu}*\varphi)^*)
|\lambda|^nd\lambda\right|^{\frac{1}{2}} = C'\|\tilde{\nu}*\varphi\|_{L^2}$$
by Plancherel formula. By the general Hausdorff-Young inequality,
$$\|\tilde{\nu}*\varphi\|_{L^2}\leq\|\tilde{\nu}\|_{L^1}\|\varphi\|_{L^2}\leq C'\|\varphi\|_{L^2}.$$
Estimation of the second integral on the right-hand side of (3.13) is similar except that (3.15) is replaced by
$$\left[\int^{\infty}_a|\lambda|^{2(N-\alpha/2)+n-N'}d\lambda\right]^{\frac{1}{2}}\left[\int^{\lambda}_0|tr(\lambda^{N'}\pi_{\lambda}(\tilde{\nu})\pi_{\lambda}
(\varphi)B_1)||\lambda|^nd\lambda\right]^{\frac{1}{2}}. \eqno (3.18)$$
Now $(i\lambda)^{N'}=\pi_{\lambda}(T^{N'})$ and $\pi_{\lambda}(T^{N'})\pi_{\lambda}(\tilde{\nu})=\pi_{\lambda}(\widetilde{T^{N'}}\nu).$\\
If $N'$ is chosen so large that the first factor is finite, (3.18) is bounded by \\
$\displaystyle C^{\prime\prime}\left\{\int^{\lambda}_0tr(\pi_{\lambda}(\widetilde{T^{N'}}\nu)\pi_{\lambda}(\varphi)\pi_{\lambda}(\varphi)^*
\pi_{\lambda}(\widetilde{T^{N'}}\nu)^*)|\lambda|^nd\lambda\right\}^{\frac{1}{2}}$
$$= C^{\prime\prime}\|\widetilde{T^{N'}}\nu*\varphi\|\leq C^{\prime\prime\prime}\|\varphi\|_{L^2}.$$

Combining these estimates and similar ones for the other integral in (3.12) gives (3.11).  This argument is very similar to that given in \cite{Rockland}.

To show that the $L^2$ function $u$ is actually smooth it suffices by Sobolev's Lemma, to show that all distributional derivatives of $u$ are in $L^2$.  We shall show that if $D$ is any left invariant differential operator of $\h_n$, then
$$|u(D\varphi)|\leq C_D\|\varphi\|_{L^2}, \;\;\varphi\in C^{\infty}_0(\h_n).\eqno (3.19)$$
To prove (3.19), one proceeds as above, but with $\varphi$ replaced by $D\varphi$. \\ Since
$tr(\pi_{\lambda}(\tilde{\nu})\pi_{\lambda}(D\varphi)\pi_{\lambda}(D\varphi)^*\pi_{\lambda}(\tilde{\nu})^*)=
tr(\pi_{\lambda}(D\varphi)^*\pi_{\lambda}(\tilde{\nu})\pi_{\lambda}(\tilde{\nu})\pi_{\lambda}(\tilde{\nu})^*
\pi_{\lambda}(D\varphi)),$\\
$\displaystyle\int_{0<\lambda<a}|\lambda|^N tr(\pi_{\lambda}(\tilde{\nu})\pi_{\lambda}(D\varphi)(\pi_{\lambda}(P))^{-1})|\lambda|^nd\lambda$
$$\leq C\|D\varphi *\tilde{\nu}\|_{L^2}= C\|\varphi *D^r\tilde{\nu}\|_{L^2}\leq C'\|\varphi\|_{L^2}.$$
Similarly,\\
$\displaystyle\int_{0\leq\lambda<\infty}|\lambda|^N tr(\pi_{\lambda}(\tilde{\nu})\pi_{\lambda}(D\varphi)(\pi_{\lambda}(P))^{-1})|\lambda|^nd\lambda$
$$\leq C\|D\varphi *T^{N'}\tilde{\nu}\|_{L^2}= C\|\varphi *D^rT^{N'}\tilde{\nu}\|_{L^2}\leq C'\|\varphi\|_{L^2}, \;\;r\in \mathbb{\R}^+.$$
The estimates for $\lambda <0$ are the same.  This completes the proof that $u$ is $C^{\infty}$ function.

Finally, a formal calculation in \cite{Rothschild1}\cite{Egwe10}\cite{BMEE1} shows that $u$ satisfies (3.8) i.e., that\\ $u(P^{\tau}\varphi)=T^N\nu(\varphi)$ for all $\varphi\in C^{\infty}_0(\h_n)$, where $P^{\tau}$ is the formal transpose of $P$.  Indeed, since
$$\pi_{\lambda}(P^{\tau}\varphi)=\pi_{\lambda}(\varphi)\pi_{\lambda}(P),$$
and
$$(i\lambda)^N\pi_{\lambda}(\tilde{\nu})=\pi_{\lambda}(\tilde{\nu})\pi_{\lambda}(T^N)=\pi_{\lambda}
(\widetilde{T^N}\nu),$$

$$\begin{array}{rcl}
u(P^{\tau}\varphi)&=&\displaystyle\int_{\mathbb{R}^*}(i\lambda)^N tr(\pi_{\lambda}(\tilde{\nu})\pi_{\lambda}(P^{\tau}\varphi)(\pi(P))^{-1})|\lambda|^Nd\lambda\\
\\
&=&\displaystyle\int_{\mathbb{R}^*} tr(\pi_{\lambda}(\widetilde{T^N}\nu)\pi_{\lambda}(\varphi))|\lambda|^Nd\lambda\\
\\
&=&\displaystyle\int_{\mathbb{R}^*} tr(T^N\nu*\varphi))|\lambda|^Nd\lambda\\
\\
&=&(T^N\nu*\varphi)|\lambda|^Nd\lambda\\
&=&(\widetilde{T^N}\nu *\varphi)(0)\;\;\;\;\;\;\mbox{(by the Plancherel theorm,)}\\
&=&T^N\nu(\varphi),\;\;\;\;\;\;;\;\;\;\;\;\;\mbox{(by (3.17))}.\end{array}$$
This completes the proof of the theorem. \hspace{0.5in} $\Box$\\

\end{document}